\documentstyle[12pt]{article}
\newtheorem{theorem}{Theorem}
\newtheorem{lemma}[theorem]{Lemma}
\newtheorem{proposition}[theorem]{Proposition}

\newtheorem{corollary}[theorem]{Corollary}
\newtheorem{remar}[theorem]{Remark}
\newenvironment{proof}{Proof:\ \ \ }{\QED}
\newenvironment{remark}{\begin{remar}\rm}{\end{remar}}
\newcommand{\QED}{{\unskip\nobreak\hfil\penalty50%
\hskip1em\hbox{}\nobreak\hfil $\Box$%
\parfillskip=0pt \finalhyphendemerits=0 \par\medskip\noindent}}
\newcommand{\bfind}[1]{\index{#1}{\bf #1}}

\newcommand{\n}{\par\noindent}
\newcommand{\sn}{\par\smallskip\noindent}

\newcommand{\bn}{\par\bigskip\noindent}

\newcommand{\parm}{\par\medskip}
\newcommand{\parb}{\par\bigskip}
\font\tenlv=msbm10 scaled 1200
\font\sevenlv=msbm7 scaled 1200
\font\fivelv=msbm5 scaled 1200
\def\lv #1{{\mathchoice{{\hbox{\tenlv #1}}}{{\hbox{\tenlv #1}}}
{{\hbox{\sevenlv #1}}}{{\hbox{\fivelv #1}}}}}
\newcommand{\N}{\lv N}

\newcommand{\Z}{\lv Z}
\newcommand{\pH}{\mathop{\rule[-1pt]{0pt}{1pt}\mbox{\large\bf H}}}
\newcommand{\subsetuneq}{\mathrel{\raisebox{.8ex}{\footnotesize%
$\displaystyle\mathop{\subset}_{\not=}$}}}

\newcommand{\bbox}[1]{\makebox(0,0){\rule[-2ex]{0ex}{5.4ex}#1}}
\font\tenlv=msbm10 scaled 1200
\font\sevenlv=msbm7 scaled 1200
\font\fivelv=msbm5 scaled 1200
\def\lv #1{{\mathchoice{{\hbox{\tenlv #1}}}{{\hbox{\tenlv #1}}}
{{\hbox{\sevenlv #1}}}{{\hbox{\fivelv #1}}}}}
\renewcommand{\N}{\lv N}
\renewcommand{\Z}{\lv Z}

\begin{document}            
\title{Functional Equations for Lexicographic Products
\footnote{2000{\it Mathematics Subject Classification}: Primary
06A05, Secondary 03C60.
\n First and second authors partially supported
by an NSERC research grant, third author partially supported by
the Edmund Landau Center for research in Mathematical Analysis,
supported by the Minerva Foundation (Germany). Publication number 615}}
\author{Franz-Viktor, Salma Kuhlmann and Saharon Shelah\\[.3cm]}
\date{3.\ 6.\ 2001}
\maketitle
%
\begin{abstract}\noindent
We generalize
the main result of [K--K--S] concerning the convex embeddings of a chain
$\Gamma$ in a lexicographic power $\Delta ^{\Gamma}$. For a fixed
non-empty chain $\Delta$, we derive
necessary and sufficient conditions for the existence of non-empty
solutions $\Gamma$ to each of the lexicographic functional equations
\[(\Delta ^{\Gamma})^{\leq 0} \simeq \Gamma\>, \ \
(\Delta ^{\Gamma}) \simeq \Gamma \ \ \mbox{ and } \ \
(\Delta ^{\Gamma})^{<0} \simeq \Gamma\;.\]
\end{abstract}

\section{Introduction}
Let us recall the definition of lexicographic products of ordered sets.
Let $\Gamma$ and $\Delta_\gamma\,$, $\gamma\in\Gamma$ be non-empty
totally ordered
sets. For every $\gamma\in\Gamma$, we fix a distinguished element
$0_\gamma\in\Delta_\gamma\,$. The \bfind{support} of a family $a=
(\delta_\gamma)_{\gamma\in\Gamma}\in\prod_{\gamma\in\Gamma}\Delta_\gamma$
is the set of all $\gamma\in\Gamma$ for which $\delta_\gamma\ne 0_\gamma$.
We denote it by $\mbox{\rm supp} (a)$. As a set, we define
$\pH_{\gamma\in\Gamma}(\Delta_\gamma , 0_\gamma)$ to be the set of all
families
$(\delta_\gamma)_{\gamma\in\Gamma}$ with well-ordered support (with
respect to fixed
distinguished elements $0_\gamma$). To relax the notation, we shall
write $\pH_{\gamma\in\Gamma}\Delta_\gamma$ instead of
$\pH_{\gamma\in\Gamma}(\Delta_\gamma , 0_\gamma)$ once the distinguished
elements $0_\gamma$ have been fixed.
Then  the
\bfind{lexicographic order} on $\pH_{\gamma\in\Gamma}\Delta_\gamma$ is
defined as follows. Given $a=(\delta_\gamma)_{\gamma\in\Gamma}$ and $b=
(\delta'_\gamma)_{\gamma\in\Gamma}\in \pH_{\gamma\in\Gamma}\Delta_\gamma$,
observe that
$\mbox{\rm supp}(a)\cup\mbox{\rm supp}(b)$ is well-ordered. Let
$\gamma_0$ be the least of all elements $\gamma\in\mbox{\rm supp}(a)\cup
\mbox{\rm supp}(b)$ for which $\delta_\gamma\ne\delta'_\gamma\,$. We set
$a<b\,:\Leftrightarrow\,\delta_{\gamma_0}<\delta'_{\gamma_0}\,$. Then
$(\pH_{\gamma\in\Gamma}\Delta_\gamma,<)$ is a totally ordered set, the
\bfind{lexicographic product} (or \bfind{Hahn product}) of the ordered
sets $\Delta_\gamma\,$.
We shall always denote by $0$ the sequence with
empty support in $\pH_{\gamma\in\Gamma}\Delta_\gamma$.

Note that if all $\Delta_\gamma$ are totally ordered
abelian groups, then we can take the distinguished elements $0_\gamma$ to be
the neutral elements of the groups $\Delta_\gamma\,$. Defining addition
on $\pH_{\gamma\in\Gamma} \Delta_\gamma$ componentwise, we obtain a
totally ordered abelian group $(\pH_{\gamma\in\Gamma}\Delta_\gamma,+,0<)$.
\sn
{\bf Lexicographic exponentiation of chains}:
If $\Delta = \Delta_\gamma$ for every $\gamma\in\Gamma$, we fix a
distinguished element in $\Delta$ (the same distinguished element for every
$\gamma\in\Gamma$), and denote it
by $0_{\Delta}$. In this case we denote $\pH_{\gamma\in\Gamma} \Delta_\gamma$
by $\Delta ^{\Gamma}$, and call it
the \bfind{lexicographic power} $\Delta ^{\Gamma}$ (with respect to
$0_{\Delta}$). In other words, $\Delta ^{\Gamma}$ is the set
\[\{s; s:\Gamma \rightarrow
\Delta \mbox{ such that\ supp($s$) is well-ordered in $\Gamma$}\},\]
ordered lexicographically.

This exponentiation of chains has its own arithmetic. In this paper we
study some of its aspects (cf.\  also [K] and [H--K--M]). Note that if
$\Gamma$ and $\Delta$ are infinite ordinals, then lexicographic
exponentiation does {\it not} coincide with ordinal exponentiation (cf.\ [H]).

Lexicographic powers appear naturally in many contexts.
For example, $\N ^{\N}$ is the order type of the nonnegative reals,
and $\Z^{\N}$ that of the irrationals (cf.\ [R]). Also, $2^{\Gamma}$
is (isomorphic to) the chain of all well-ordered subsets of $\Gamma$,
ordered by inclusion. The chain $2^{\N}$ has been studied in [H].

However, the main motivating example for us was that of generalized
power series fields. If $k$ is a real closed field and $G$ a totally
ordered divisible abelian group, then the field $k((G))$ of power series
with exponents in $G$ and coefficients in $k$ is again real closed. The
unique order of $k((G))$ is precisely the chain $k^G$. It was while
studying such fields that our interest in the present problems arose. In
[K-K-S], we considered the problem of defining an exponential function
on $K=k((G))$, that is, an isomorphism $f$ of ordered groups $f: (K, +,
0, <) \rightarrow (K^{>0}, \cdot, 1, <)$. We showed that the existence
of $f$ would imply that of a \bfind{convex embedding} (that is, an
embedding with convex image) of the chain $G^{<0}$ into the chain
$k^{G^{<0}}$. On the other hand, we proved:
\begin{theorem}                              \label{previous}
Let $\Gamma$ and $\Delta$ be non-empty totally ordered sets without
greatest element, and fix an element $0_\Delta\in\Delta$. Suppose that
$\Gamma'$ is a cofinal subset of $\Gamma$ and that $\iota\colon
\Gamma'\,\rightarrow\,\Delta^\Gamma$ is an order preserving embedding.
Then the image $\iota\Gamma'$ is not convex in $\Delta^\Gamma$.
\end{theorem}
Now for any ordered field $k$, the chain $k$ has no last element.
Similarly, $G^{<0}$ has no last element if $G$ is nontrivial and
divisible. So, using Theorem~\ref{previous} one establishes that no
exponentiation is possible on generalized power series fields.

If we omit the conditions on $\Gamma$ and $\Delta$ in
Theorem~\ref{previous}, the situation changes drastically. In this
paper, we study conditions on the chains $\Gamma$ and $\Delta$ under
which a convex embedding of $\Gamma$ in $\Delta ^{\Gamma}$ exists. In
particular, we seek for non-empty solutions $\Gamma$ to the functional
equations:
\[(\Delta ^{\Gamma})^{\leq 0} \simeq \Gamma\,, \
(\Delta ^{\Gamma}) \simeq \Gamma\,, \ \mbox{ and } \
(\Delta ^{\Gamma})^{<0} \simeq \Gamma\;.\]
(if $T$ is any totally ordered set and $0\in T$ is any element, we
denote by $T^{\leq 0}$ the initial segment (including $0$),
and by $T^{<0}$ the strict initial segment (excluding $0$) determined by
$0$ in $T$).
None of the three equations hold if both $\Delta$ and $\Gamma$ have no
last element (for the first, this is trivial, and for the second and
third it follows from Theorem~\ref{previous}). In
Section~\ref{sectgen} we start by proving a strong generalization of
Theorem~\ref{previous} (cf.\ Theorem \ref{2001}).
In Section~\ref{sectsol}, for each of the three functional equations,
we give simple characterizations of those
chains $\Delta$ for which non-empty solutions ${\Gamma}$ exist.
In Section~\ref{simul} we study simultaneous solutions to all three
equations.

%
%
\section{Nonexistence of convex embeddings}          \label{sectgen}
In this section, we shall prove that Theorem~\ref{previous} remains true
in the case where $\Delta$ is arbitrary, but $0_\Delta$ is not the last
element of $\Delta$. This will follow from the following more general
result:
\begin{theorem}                           \label{2001}
Let $\Gamma$ and $\Delta_\gamma\,$, $\gamma\in\Gamma\,$, be
non-empty totally
ordered sets. For every $\gamma\in\Gamma\,$, fix an element $0_\gamma$
which is not the last element in $\Delta_\gamma\,$. Suppose that
$\Gamma$ has no last element and that $\Gamma'$ is a cofinal subset of
$\Gamma$. Then there is no convex embedding
\[\iota:\>\Gamma'\,\rightarrow\,\pH_{\gamma\in\Gamma} \Delta_\gamma\:.\]
\end{theorem}
\begin{proof}
For every $\gamma\in\Gamma'$, we choose an element $1_\gamma\in
\Delta_\gamma$ such that $1_\gamma>0_\gamma\,$. Take $d=
(d_\gamma)_{\gamma\in\Gamma}\,$. If $S$ is a well-ordered subset of
$\Gamma'$ such that $d_\gamma=0_\gamma$ for all $\gamma\in S$, then we
set
\[d\oplus S\>:=\> (d'_\gamma)_{\gamma\in\Gamma}\mbox{ \ \ with \ \ }
d'_\gamma = \left\{\begin{array}{rl}
d_\gamma & \mbox{for } \gamma\notin S\\
1_\gamma & \mbox{for } \gamma\in S\;.
\end{array}\right.\]
Observe that the support of $d\oplus S$ is contained in
$\mbox{\rm supp}(d)\cup S$ and thus, it is again well-ordered.
Note also that
\begin{equation}                            \label{oplusS}
S'\subsetuneq S\>\Rightarrow\> d\oplus S' < d\oplus S\;.
\end{equation}
Indeed, let $\gamma_0$ be the least element in $S\setminus S'$. Then
$(d\oplus S')_{\gamma _0} = 0_\gamma < 1_\gamma = (d\oplus
S)_{\gamma_0}$. On the other hand, if $\gamma \in \Gamma$ and $\gamma <
\gamma _0$ then $(d\oplus S')_\gamma = (d\oplus S)_\gamma$: if $\gamma
\in S$ then $\gamma \in S'$ (by minimality of $\gamma _0$) and
$(d\oplus S')_\gamma = 1_\gamma = (d\oplus S)_\gamma$; if
$\gamma \notin S$ then $\gamma \notin S'$ and
$(d\oplus S')_\gamma = d_\gamma = (d\oplus S)_\gamma$.

Now suppose that
$\iota:\>\Gamma'\,\rightarrow\,\pH_{\gamma\in\Gamma} \Delta_\gamma$
is an order preserving embedding such that the image
$\iota\Gamma'$ is convex in $\pH_{\gamma\in\Gamma} \Delta_\gamma$. We
wish to deduce a
contradiction. The idea of the proof is the following. Let ON denote the
class of ordinal numbers. We shall define an infinite $\mbox{\rm ON}
\times\N$ matrix with coefficients in $\Gamma'$, such that each
column $(\gamma_{\nu}^{(n)})_{\nu\in {\rm ON}}$ is a strictly
increasing sequence in $\Gamma'$. Since $\Gamma'$ is a {\it set}, every
column of this matrix will provide a contradiction at the end of the
construction (cf.\ figure).

\begin{center}
\setlength{\unitlength}{0.002\textwidth}
\begin{picture}(200,200)(-100,-100)
\put(0,0){\bbox{$\left(\begin{array}{ccccc}
\gamma_0^{(1)} & \ldots & \gamma_0^{(n)} & \gamma_0^{(n+1)}\; &
\;\ldots \\
\vdots & & \vdots & \vdots & \\[-7pt]
\vdots & & \vdots & \vdots & \\
\gamma_{\nu}^{(1)} & \ldots & \gamma_{\nu}^{(n)} & \gamma_{\nu}^{(n+1)}
& \ldots \\
\vdots & & \vdots & \vdots & \\
\gamma_{\mu}^{(1)} & \ldots & \gamma_{\mu}^{(n)}=\,? & \ldots & \ldots
\\
\vdots & & & & \\
\ldots & \ldots & \ldots & \ldots & \ldots
\end{array}\right)$}}
\put(-30,62){\line(1,0){50}}
\put(-30,92){\line(1,0){105}}
\put(-30,62){\line(0,1){30}}
\put(20,62){\line(0,-1){87}}
\put(75,92){\line(0,-1){117}}
\put(20,-25){\line(1,0){55}}
\end{picture}
\end{center}
To get started, we have to define the first row of the matrix.
We construct sequences $\beta^{(n)}$, $n\in\N\cup\{0\}$, and
$\gamma_0^{(n)}$, $n\in\N$, in $\Gamma'$. We take an arbitrary
$\beta^{(0)}\in\Gamma'$. Having constructed $\beta^{(n)}$, we choose
$\gamma_0^{(n+1)}$ and $\beta^{(n+1)}$ as follows. Since $\Gamma'$ has
no last element, we can choose $\mu^{(n)},\nu^{(n)}\in\Gamma'$ such that
$\beta^{(n)}<\mu^{(n)}<\nu^{(n)}$. Hence,
\[\iota\beta^{(n)}<\iota\mu^{(n)}<\iota\nu^{(n)}\;.\]
Let $\sigma^{(n)}\in\Gamma$ be the least element in $\mbox{\rm supp}
\,\iota\beta^{(n)}\cup\mbox{\rm supp}\,\iota\mu^{(n)}$ for which
\begin{equation}                            \label{sig}
(\iota\beta^{(n)})_{\sigma^{(n)}}<(\iota\mu^{(n)})_{\sigma^{(n)}}\;,
\end{equation}
and $\tau^{(n)}\in\Gamma$ the least element in $\mbox{\rm supp}
\,\iota\mu^{(n)}\cup\mbox{\rm supp}\,\iota\nu^{(n)}$ for which
\begin{equation}                            \label{tau}
(\iota\mu^{(n)})_{\tau^{(n)}}<(\iota\nu^{(n)})_{\tau^{(n)}}\;.
\end{equation}
Since $\Gamma'$ is cofinal in $\Gamma$, we can choose $\beta^{(n+1)}
\in\Gamma'$ such that
\[\beta^{(n+1)}\,\geq\,\max\{\sigma^{(n)},\tau^{(n)}\}\;.\]
Further, we set
\[d^{(n+1)}\>:=\>(d^{(n+1)}_\gamma)_{\gamma\in\Gamma}\mbox{ \ \ with \ \ }
d^{(n+1)}_\gamma = \left\{\begin{array}{rl}
(\iota\mu^{(n)})_\gamma & \mbox{for } \gamma\leq\beta^{(n+1)}\\
0_\gamma & \mbox{for } \gamma>\beta^{(n+1)}\;.
\end{array}\right.\]
Then by (\ref{sig}) and (\ref{tau}),
\[\iota\beta^{(n)}<d^{(n+1)}<\iota\nu^{(n)}\;.\]
Thus, $d^{(n+1)}\in\iota\Gamma'$ by convexity, and we can set
\[\gamma_0^{(n+1)}\>:=\>\iota^{-1}d^{(n+1)}\;.\]
Now for every $n\in\N$ we have that $\beta^{(n)}<\gamma_0^{(n+1)}\,$,
hence every well-ordered set $S\subset\Gamma'$ with smallest element
$\gamma_0^{(n+1)}$ has the property that $(\iota\gamma_0^{(n)})_\gamma=
d^{(n)}_\gamma=0_\gamma$ for all $\gamma\in S$; and moreover,
\[\iota\gamma_0^{(n)}<\iota\gamma_0^{(n)}\oplus S<\iota\nu^{(n-1)}\;.\]
Thus, $\iota\gamma_0^{(n)}\oplus S\in\iota\Gamma'$ by convexity.
Suppose now that for some ordinal number $\mu\geq 1$ we have chosen
elements $\gamma_\nu^{(n)}\in\Gamma'$, $\nu<\mu$, $n\in\N$, such that
for every fixed $n$, the sequence $(\gamma_\nu^{(n)})_{\nu<\mu}$ is
strictly increasing. Then we set
\[\gamma_\mu^{(n)}:=\iota^{-1}(\iota\gamma_0^{(n)}\oplus
\{\gamma_\nu^{(n+1)} \mid\nu<\mu\})\in\Gamma'\]
for every $n\in\N$. If $\lambda<\mu$, then $\{\gamma_\nu^{(n+1)} \mid
\nu<\lambda\}\subsetuneq\{\gamma_\nu^{(n+1)} \mid\nu<\mu\}$ and thus,
$\gamma_\lambda^{(n)}<\gamma_\mu^{(n)}$ by (\ref{oplusS}).
So for every ordinal number $\mu$, the sequences
$(\gamma_\nu^{(n)})_{\nu<\mu}$ can be extended. We obtain strictly
increasing sequences of arbitrary length, contradicting the fact that
their length is bounded by the cardinality of $\Gamma$.
\end{proof}

\begin{corollary}
Assume that $0_\Delta$ is not the last element of $\Delta$. If there is
an embedding of $\Gamma$ in $\Delta^\Gamma$ with convex image, then
$\Gamma$ has a last element.
\end{corollary}
%
%
\section{Solutions to the Functional equations}    \label{sectsol}
We start with a few easy remarks and lemmas.
Throughout, fix a chain $\Delta$ with distinguished element
$0_{\Delta}$.
\begin{remark}               \label{r1}
1) If $0_{\Delta}$ is last in $\Delta$
(respectively, least), then $0$ is last in $\Delta ^{\Gamma}$
(respectively, least), for any non-empty chain $\Gamma$.\n
2) Let $I$ be any chain, and $C$ a non-empty convex subset of $I$.
Let $c\in C$. Then the initial segment determined by $c$ in $C$ is a
final segment of the initial segment determined by $c$ in $I$.
\end{remark}
\begin{remark}                    \label{r2}
If  $\Delta ^{<0_\Delta}$ has no last element, then also
$(\Delta ^{\Gamma}) ^{<0}$ has no last, for any chain $\Gamma$:
If not, let $s$ be last in $(\Delta ^{\Gamma}) ^{<0}$
and set $\gamma = \mbox{min supp} (s)$.
Then $s(\gamma) = \delta < 0_{\Delta}$.
Take $\delta < \delta ' <0_{\Delta}$. Consider $s'$ defined by
$s'(\gamma)= \delta '$ and $s'(\gamma ')=0_{\Delta}$ if $\gamma '\not=
\gamma$. Then $s'\in (\Delta ^{\Gamma}) ^{<0}$, but $s' > s$,
contradiction.\n
\end{remark}

\begin{lemma}                          \label{lifting}
Let $\Gamma$ and $\Gamma '$ be chains, and suppose that
$\phi : \Gamma \rightarrow \Gamma '$ is a chain embedding.
Then $\phi$ lifts to a chain embedding
\[\hat{\phi} :\Delta ^{\Gamma} \rightarrow \Delta ^{\Gamma '}\;. \]
\end{lemma}
\begin{proof}
For $s \in \Delta ^{\Gamma}$ and $x \in \Gamma '$, set
\[ \hat {\phi} (s) (x) = \left\{\begin{array}{ll}
0_{\Delta} & \mbox{if } x\notin \mbox{Im}\,\phi\\
s(\phi ^{-1}(x)) & \mbox{if } x \in \mbox{Im}\,\phi\;.
\end{array}\right.\]
(here, Im $\phi$ denotes the image of $\phi$).
Now, it is straightforward to check the assertion of the lemma.
\end{proof}

In view of this lemma, if $F$ is a subchain of a chain $\Gamma$,
then there is a natural identification of $\Delta ^F$ as a subchain of
$\Delta ^{\Gamma}$.
\begin{lemma}                                \label{final segment}
Let $\Gamma$ be a chain and $F$ a non-empty final segment of $\Gamma$.
Then $\Delta ^F$ is convex in $\Delta ^{\Gamma}$ (and $0\in
\Delta ^F$).
\end{lemma}
\begin{proof}
Let $s_i \in \Delta ^F$, and set $\gamma _i= \mbox {min supp} (s_i) \in
F$, for $i=1,2$. Let $s\in \Delta ^{\Gamma}$ be such that $s_1 < s <
s_2$. If $s=0$, then $s\in \Delta ^F$. So assume $s\not= 0$ and set
$\gamma = \mbox {min supp} (s)$. Suppose that $\gamma \notin F$. If
$s>0$, then $s(\gamma) >0_{\Delta}$. On the other hand, $\gamma < \gamma
_2$ (otherwise, $\gamma \in F$). Thus, $s > s_2\,$, a contradiction.
Similarly, we argue that if $s<0$, then $s < s_1\,$, a contradiction.
Hence, $\mbox {min supp} (s)$. Since $F$ is a final segment of $\Gamma$,
this implies that $s\in \Delta ^F$, which proves our assertion.
\end{proof}

\begin{corollary}           \label{last}
Assume that $\Gamma$ has a last element. Then $\Delta$ embeds convexly
in $\Delta ^{\Gamma}$, such that $0_{\Delta}$ is mapped to $0 \in
\Delta ^{\Gamma}$. If moreover $0_{\Delta}$ is last in $\Delta$, then
$\Delta ^F$
embeds as a final segment in $\Delta ^{\Gamma}$, for any non-empty final
segment
$F$ of $\Gamma$. Consequently, if $\Gamma$ has a last element, and
$0_{\Delta}$ is last in $\Delta$, then $\Delta$ embeds as a final
segment in $\Delta ^{\Gamma}$.
\end{corollary}
\begin{proof}
The first assertion follows from Lemma \ref{final segment}, applied to
the final segment consisting of the single last element of $\Gamma$.
For the second assertion use
Remark \ref{r1}, parts 1) and 2).
\end{proof}

\parm
We now give a complete solution to the {\bf first functional equation},
and a sufficient condition for the existence of solutions $\Gamma$ to
the third functional equation:
\begin{theorem}                            \label{equation1}
There is always a non-empty solution $\Gamma$ for the functional equation
$(\Delta^{\Gamma})^{\leq 0} \simeq \Gamma$.
If $\Delta^{<0_{\Delta}}$ has a last element, then there is also a non-empty
solution $\Gamma$ for $(\Delta^{\Gamma})^{<0} \simeq \Gamma$.
\end{theorem}
\begin{proof}
Set $\Gamma_0 := \Delta^{\leq 0_{\Delta}}$. Since $\Gamma_0$ has a last
element, $\Delta$ embeds convexly in $\Delta^{\Gamma_0}$. Consequently,
$\Gamma_0$ embeds as a final segment in $\Gamma_1:=(\Delta^{\Gamma
_0})^{\leq 0}$. By induction on $n\in \N$ we define $\Gamma_n:=
(\Delta^{\Gamma _{n-1}})^{\leq 0}$, and obtain an embedding of $\Gamma
_{n-1}$ as a final segment in $\Gamma _n$. We set $\Gamma :=
\cup _{n\in \N} \Gamma _n$.

Since every $\Gamma_n$ is a final segment of $\Gamma$, every
well-ordered subset $S$ of $\Gamma$ is already contained in some
$\Gamma_n$ (just take $n$ such that the first element of $S$ lies in
$\Gamma_n$). Hence, an element of $(\Delta^\Gamma)^{\leq 0}$ with
support $S$ is actually an element of
$\Gamma_{n+1}=(\Delta^{\Gamma_n})^{\leq 0}$, for some $n$. This fact
gives rise to an order isomorphism of $(\Delta^\Gamma)^{\leq 0}$ onto
$\Gamma$.

To prove the second assertion, we set $\Gamma_0:=\Delta^{<0_{\Delta}}$.
Since $\Gamma_0$ has a last element by assumption, $\Delta$ embeds
convexly in $\Delta^{\Gamma_0}$, and the same arguments as above work
if we define $\Gamma_n:= (\Delta^{\Gamma _{n-1}})^{< 0}$.
\end{proof}
\begin{remark}
Note that $\Gamma_0$ has a last element and embeds as a final segment in
the constructed solution $\Gamma$ (in both cases considered in the
proof). Thus, $\Gamma$ has a last element, and there is no contradiction
to Theorem \ref{2001}.
%
\end{remark}

Note that if $0_{\Delta}$ is least in $\Delta$, then the first equation
has the trivial solution $\Gamma = \{0_{\Delta}\}$.

\parb
We next turn to the {\bf second functional equation}.
\begin{remark}                        \label{Mai}
Suppose that $0_{\Delta}$ is last in $\Delta$. Then the solution to the
first equation given in Theorem \ref{equation1} also solves the second
equation. Indeed, in this case, $0$ is last in $\Delta ^{\Gamma}$, so
$(\Delta ^{\Gamma})^{\leq 0}=\Delta ^{\Gamma}$.
\end{remark}

We also have the converse:
\begin{corollary}                      \label{equation2}
Assume $\Delta$ is a chain such that the functional equation
$\Delta ^{\Gamma} \simeq \Gamma$ has a non-empty solution $\Gamma$. Then
$0_{\Delta}$ is last in $\Delta$. Thus, the functional equation $\Delta
^{\Gamma} \simeq \Gamma$ has a non-empty solution if and only if
$0_{\Delta}$ is last in $\Delta$.
\end{corollary}
\begin{proof}
Assume $0_{\Delta}$ is not last, and choose some element $1_{\Delta}>
0_{\Delta}$. This provides us with characteristic functions. If $S
\subset \Gamma$ is well-ordered, then let $\chi _S\in \Delta ^{\Gamma}$
denote the characteristic function on $S$ defined by:
\[\chi _S (\gamma) =\left\{\begin{array}{ll}
1_\Delta & \mbox{if } \gamma \in S\\
0_{\Delta} & \mbox{if } \gamma \notin S\;.
\end{array}\right.\]
Note that these characteristic functions reflect inclusion: if $S$ is a
proper well-ordered subset of $S'$, then $\chi _S < \chi_{S'}$.
Now assume for a contradiction that $i: \Gamma\simeq \Delta ^{\Gamma}$,
and let $\kappa = \mbox{card} (\Gamma)$. We shall construct a strictly
increasing sequence $\{\gamma _{\mu}; \mu < \kappa ^+\}$ in~$\Gamma$.

Set $\gamma _0 = i^{-1}(0)$, and assume by induction that
$\{\gamma _{\nu}; \nu < \mu \}$ is defined, and strictly increasing in
$\Gamma$.
Then define
\[\gamma _{\mu}= i^{-1}(\chi_ {\{\gamma _{\nu}; \nu < \mu \}}).\]
It follows that $\chi_ {\{\gamma _{\lambda}; \lambda < \nu \}}
< \chi_ {\{\gamma _{\lambda}; \lambda < \mu \}}$, whenever $\nu < \mu$.
Since $i^{-1}$ is order preserving,
it follows that $\gamma _{\nu} < \gamma _{\mu}$ as required.
\end{proof}

\parm
We now turn to the {\bf third functional equation}.
We deduce a simple criterion for the existence of solutions:
\begin{corollary}                          \label{11.Mai}
Assume that $0_\Delta$ is not the last element of $\Delta$.
Then the functional equation $(\Delta^{\Gamma})^{<0} \simeq \Gamma$ has
a non-empty solution $\Gamma$ if and only if $\Delta^{<0_{\Delta}}$ has a
last element.
\end{corollary}
\begin{proof}
The ``if'' direction is just the second assertion of Theorem
\ref{equation1}. So assume now that $\Gamma$ is a non-empty solution.
Assume for a contradiction that
$\Delta^{<0_{\Delta}}$ has no last element. Then by Remark
\ref{r2} $(\Delta^{\Gamma})^{<0}$ has no last element as well. Thus,
the same holds for the solution $\Gamma$. This contradicts Theorem
\ref{2001}.
\end{proof}

%
%
\section{Simultaneous Solutions}          \label{simul}
Recall that by Remark~\ref{Mai}, the chain $\Gamma$ given in Theorem
\ref{equation1} solves the first {\it and} the second functional
equations, if $0_{\Delta}$ is last in $\Delta$.

\begin{theorem}                          \label{sufficient}
Assume that $0_{\Delta}$ is last in $\Delta$ and that $\omega ^*$
embeds as a final segment in $\Delta$. Then the solution $\Gamma$
given in Theorem \ref{equation1} to the first and second functional
equations solves
$(\Delta ^{\Gamma}) ^{<0} \simeq \Gamma$ as well.
\end{theorem}
\begin{proof}
Recall that $\Delta$ embeds as a final segment in the given solution
$\Gamma$. Thus, $\omega ^*$ embeds as a final segment in $\Gamma$ as
well. In particular, $\Gamma$ has a last element $0$. Since
$\Delta ^{\Gamma} = (\Delta ^{\Gamma}) ^{<0} \cup \{0\}$ and
$\Delta ^{\Gamma} \simeq \Gamma$, we find that
$(\Delta ^{\Gamma})^{<0} \simeq \Gamma \setminus \{0\}$. But $\Gamma
\simeq \Gamma \setminus \{0\}$, since $\omega ^*$ is a final segment of
$\Gamma$.
\end{proof}

\parm
We now turn to the question of whether the sufficient conditions given
in this last theorem is also necessary. We need to introduce a
definition: Say that a solution $\Gamma$ (to any of the three equations)
is \bfind{special} if $\Delta$ embeds as a final segment in $\Gamma$.
Note that special solutions are necessarily non-empty.

\begin{proposition}                           \label{special}
Every non-empty solution to $\Gamma \simeq \Delta ^{\Gamma}$ is special.
\end{proposition}
\begin{proof}
Necessarily, $0_{\Delta}$ is last in $\Delta$ (by
Corollary~\ref{equation2}). Thus, $\Gamma$ has a last element, so by
Corollary~\ref{last}, $\Delta$ embeds as a final segment in $\Delta
^{\Gamma}$, and thus in $\Gamma$.
\end{proof}

\begin{corollary}               \label{partialanswer}
Assume that $\Delta$ is infinite and $\Gamma$ is any non-empty chain
which solves simultaneously
\[(\Delta ^{\Gamma})^{<0} \simeq \Gamma \simeq \Delta ^{\Gamma}.\]
Then $0_{\Delta}$ is last in $\Delta$ and $\omega ^*$ embeds as a final
segment in $\Delta$.
\end{corollary}
\begin{proof}
Since $\Gamma \simeq \Delta ^{\Gamma}$, $0_{\Delta}$ is last in $\Delta$
(Corollary \ref{equation2}). Therefore, $0$ is last in $\Delta^{\Gamma}$
by Remark \ref{r1}, and so also $\Gamma$ has a last element $0$. The
assumptions imply that $\Gamma \setminus \{0\} \simeq \Gamma$.
This is equivalent to the assertion that $\omega ^*$ embeds as a final
segment in $\Gamma$. Now note that $\Gamma$ is a special solution by
Proposition \ref{special}, i.e., $\Delta$ embeds as a final segment of
$\Gamma$. Since $\Delta$ is infinite this implies that $\omega ^*$
embeds as a final segment in $\Delta$, as required.
\end{proof}

\begin{corollary}              \label{beersheva}
Assume that $\Delta$ is infinite. Then the following are
equivalent:
\n
(a) \ $0_{\Delta}$ is last in $\Delta$ and $\omega ^*$ embeds as a final
segment in $\Delta$.
\n
(b) \ There exists a (special) simultaneous solution to all three
equations.
\n
(c) \ There exists a (special) simultaneous solution to the second and
third equations.
\end{corollary}
\begin{proof}
(a) implies (b) by Theorem \ref{sufficient}. (b) implies (c) trivially.
Finally, (c) implies (a) by Corollary \ref{partialanswer}.
\end{proof}

We conclude with the following question:
%
{\it Are special solutions unique up to isomorphism?}
We can give a partial answer to this last question:
\begin{proposition}
Assume that $0_{\Delta}$ is last in $\Delta$.
Let $\Gamma=\cup \Gamma _n$ be the solution to the second equation given
in Theorem \ref{equation1}. Then $\Gamma$ embeds as a final segment in
any other solution.
\end{proposition}
\begin{proof}
Let $\Gamma '$ be another solution. Then it is a special solution, by
Proposition \ref{special}. So $\Delta =\Gamma _0$ embeds as a final
segment in $\Gamma '$. Since $0_{\Delta}$ is last in $\Delta$, $\Gamma
_1=\Delta ^{\Gamma_0}$ embeds as a final segment in $\Delta^{\Gamma'}$.
By induction, $\Gamma _n$ is a final segment of $\Gamma '$ for every
$n\in \N$. Thus. $\Gamma$ embeds as a final segment in $\Gamma '$ as
well.
\end{proof}
\bn
{\bf References}
{\small\rm
\newenvironment{reference}%
{\begin{list}{}{\setlength{\labelwidth}{5em}\setlength{\labelsep}{0em}%
\setlength{\leftmargin}{5em}\setlength{\itemsep}{-1pt}%
\setlength{\baselineskip}{3pt}}}%
{\end{list}}
\newcommand{\lit}[1]{\item[{#1}\hfill]}
\begin{reference}
\lit{[H]} {Hausdorff, F.$\,$: Grundz\"uge der Mengenlehre, Verlag von
Veit, Leipzig (1914)}
\lit{[H--K--M]} {Holland, W.~C.\ -- Kuhlmann, S.\ -- McCleary, S.$\,$:
The Arithmetic of Lexicographic Exponentiation, preprint}
\lit{[K]} {Kuhlmann, S.$\,$: Isomorphisms of Lexicographic Powers of the
Reals, Proc.\ Amer.\ Math.\ Soc.\ {\bf 123} (1995), 2657-2662}
\lit{[K--K--S]} {Kuhlmann, F.-V.\ -- Kuhlmann, S.\ -- Shelah, S.$\,$:
Exponentiation in power series fields,
Proc.\ Amer.\ Math.\ Soc.\ {\bf 125} (1997), 3177-3183}
\lit{[R]} {Rosenstein, J.~G.$\,$: Linear orderings, Academic Press,
New York - London (1982)}
\end{reference}}

\bn
\bn
\parbox[t]{8cm}{\small\rm
Department of Mathematics and Statistics\n
University of Saskatchewan\n
106 Wiggins Road\n
Saskatoon, SK S7N 5E6}
\hfil
\parbox[t]{6.5cm}{\small\rm
Department of Mathematics\n
The Hebrew University of Jerusalem\n
Jerusalem, Israel}
\end{document}